\documentclass[10pt]{article}

\usepackage{amssymb,amsthm,amsmath}

\theoremstyle{plain}
\newtheorem{theorem}{Theorem}[section]
\newtheorem{lemma}[theorem]{Lemma}
\newtheorem{corollary}[theorem]{Corollary}

\theoremstyle{definition}
\newtheorem{definition}[theorem]{Definition}
\theoremstyle{remark}

\newcommand{\auf}{\left\langle}
\newcommand{\zu}{\right\rangle}

\def\proof{\trivlist \item[\hskip \labelsep{\bf Proof:\ }]}
\def\endproof{\null\hfill\qed\endtrivlist}

\def\sgn{{\rm sgn}}
\def\const{{\rm const}}

\newcommand{\ol}{\overline}
\newcommand{\pa}{\partial}
\newcommand{\ve}{\varepsilon}
\newcommand{\R}{{\mathbb R}}

\newcommand{\cD}{\cal D}

\def\O{{\mathcal O}}

\def\R{I\!\!R}

\def\eps{\varepsilon}
\def\od{\stackrel{\rm def}{=}}

\begin{document}

\title{Weak asymptotic of shock wave formation process}
\author{V.~Danilov\footnote{Moscow Technical University of
Communication and Informatics, Russia,\hfill\break
danilov@miem.edu.ru} \and D.~Mitrovic\thanks{Faculty of
Mathematics and Natural Sciences, University of Montenegro,
matematika@cg.yu}}

\date{}

\maketitle

\begin{abstract}
We construct an asymptotic (in a weak sense) solution
corresponding to the shock wave formation in a special situation.
\end{abstract}

\section{Introduction}

We consider the problem of shock wave formation for the following
Hopf type equation:
\begin{equation}
\frac{\pa u}{\pa t}+\frac{\pa }{\pa x}f(u)=0,
\end{equation}
where we assume that $f\in C^3$ and the inequality $f''(u)>0$
holds on the range of the solution~$u$. We shall consider the
special initial condition for Eq.~(1):
\begin{equation}
u\big|_{t=0}=u^0_0+(u_1(x)-u^0_0)H(a_1-x)+(U-u_1(x))H(a_2-x),
\end{equation}
where $u^0_0$, $U$, and $a_1>a_2$ are constants, $H$ is the
Heaviside function, the function $u_1(x)$ is determined by the
equation
\begin{equation}
f'(u_1(x))=-Kx+b,\qquad K,b=\const,
\end{equation}
and, in addition, we assume that $u_1(a_1)=u^0_0$ and
$u_1(a_2)=U$.

Such a function appears\footnote{E.~Yu.~Panov drew the author's
attention to this fact.} in the construction of the entropy
solution to the Cauchy problem with an "unstable" initial jump.

It follows from the choice of such an initial condition that the
approximation of problem (1)--(2) (a weak asymptotic solution) for
all~$t$ is an element of the asymptotic subalgebra
$$
{\cal B}\{1,H_1(x-\varphi_1,\ve),H_2(x-\varphi_2,\ve)\}
$$
introduced in~\cite{1}.

Roughly speaking, this means that at any time moment the weak
limit of the weak asymptotic solution is a linear combination 
of the Heaviside functions $H(x-\varphi_1)$ and
$H(x-\varphi_2)$ with smooth in~$t$ ($\ve>0$) coefficients and
that there are no additional jumps. In turn, this means that at
time
$$
t^*=\frac{a_1-a_2}{f'(U)-f'(u^0_0)}
$$
all characteristics meet at the same point $x^*=a_i+V_i t^*$,
$V_1=f'(u^0_0)$, $V_2=f'(U)$, $i=1,2$.

More precisely, for $0<t<t^*$, the solution of problem (1)--(2) is
given by the formula
\begin{equation}
u=u^0_0+\big(u_1 (x_0(x,t))-u^0_0\big)H(\varphi_1-x)
+\big(U-u_1(x_0(x,t))\big)H(\varphi_2-x),
\end{equation}
where the function $u_1(x_0(x,t))$ has the form
$$
u_1(x_0(x,t))=u_1\bigg(a_i+\frac{x-\varphi_i(t)}{\psi_0}\psi^0_0\bigg)
=u_1\bigg(\frac{x-bt}{1-Kt}\bigg).
$$
Here $\varphi_i(t)=a_i+V_it$, $i=1,2$,
$\psi_0=\varphi_1(t)-\varphi_2(t)$, and $\psi^0_0=a_1-a_2$.

For $t=t^*$ the plot of the function $u=u(x,t)$ is the graph
$$
((-\infty,x^*),U)\cup(x^*,(U,u^0_0))\cup((x^*,\infty),u^0_0).
$$
We note that if we set
$$
\ol{u}(x,t)=\begin{cases}
u_1(x_0(x,t)),&t<t^*,\quad t>t^*,\\
u^0_0,        &t=t^*,\quad x<x^*,\\
U,            &t=t^*,\quad x>x^*,\\
\ol{u}\in[u^0_0,U],&t=t^*,\quad x=x^*,
\end{cases}
$$
then the function $\ol{u}(x,t)$ is defined  for all values of~$t$
and, for $t<t^*$, is a solution of Eq.~(1) satisfying the initial
condition~(2) for $t=0$. Our goal is to "correct" the function
$\ol{u}(x,t)$ and to obtain an analytic formula that, for $t<t^*$,
determines a function close to $\ol{u}(x,t)$ and, for $t>t^*$, a
function close to the function
\begin{equation}
u= u^0_0+H(c(t-t^*)-(x-x^*))U,
\end{equation}
where
$$
c=\frac{[f(u)]}{[u]}\bigg|_{x=ct} =\frac{f(U)-f(u^0_0)}{U-u^0_0}.
$$

The answer is given by formula (10) below.

We note that the function $u$ determined by relation~(4) for
$t<t^*$ is continuous everywhere except the points lying on the
curves $x=\varphi_i(t)$, $i=1,2$, and, at points of these curves,
the function has weak discontinuities (the derivatives of the
function have jumps at these points). Therefore, the formation of
the shock wave~(5) from~(4) can be treated as the result of
interaction (confluence) of weak discontinuities. Moreover, for
$t<t^*$, although the derivatives are discontinuous, the solution
of such problems (that is continuous, but with jumps of the
derivatives on some smooth nonintersecting curves) can be
constructed by the method of characteristics.

We note that, in the case $f(u)=u^2$, the problem of constructing
the global asymptotic solution of problem~(1), (2) was solved as
an example in~\cite{2}. The asymptotic solution constructed
in~\cite{2} is a weak asymptotic solution. We recall how it is
determined. By $O_{\cD'}(\ve^\alpha)$ we denote generalized
functions that, in general, depend on the parameters~$t$ and~$\ve$
and are such that for any test function $\eta(x)$, the estimate
$$
\auf O_{\cD'}(\ve^\alpha),\eta(x)\zu=O(\ve^\alpha)
$$
holds, where the estimate on the right-hand side is understood in
the usual sense and locally uniform in~$t$, i.e.,
$|O(\ve^\alpha)|\leq C_T\ve^\alpha$ for $t\in[0,T]$.

A function $u_\ve=u_\ve(x,t)$ is called a weak asymptotic solution
of problem~(1), (2) if
\begin{equation}
\frac{\pa u_\ve}{\pa t}+\frac{\pa f(u_\ve)}{\pa x}=O_{\cD'}(\ve),
\qquad u_\ve\bigg|_{t=0}-u\bigg|_{t=0}=O_{\cD'}(\ve).
\end{equation}

The goal of this paper is to construct such a function in the case
of a general convex nonlinearity $f(u)$. This is achieved in
Sections~2 and~3.

In Section 4 we introduce auxiliary formulas and statements of the
weak asymptotic method.

We note that if the solution $u_\ve$ satisfies the
Oleinik--Kruzhkov stability conditions~\cite{3,4}, then it follows
from~(6) that~$u_\ve$ differs from~$u$ by a measure~\cite{5} whose
values are estimated as $O(\ve)$. Indeed, it is easy to verify
that the right-hand sides in~(6) arising in our construction
belong to $C([0,T],L^1(\R^1_x))$ and can be estimated as $O(\ve)$
in the sense of the $L^1$-norm. Therefore, according to the
results in~\cite{3,4}, $u_\ve$ is an asymptotic of the solution to
the Cauchy problem (1)--(2) in~$L^1$. This is done in Section~5.

We also note that the asymptotic (in the usual weak sense) solution
describing the global behavior of the solution of the Cauchy
problem with a small viscosity and a smooth initial condition for
the equation
$$
\frac{\pa u}{\pa t}+\frac{\pa f(u)}{\pa x} =\ve\frac{\pa^2 u}{\pa
x^2}
$$
was first constructed by A.~M.~Il'in~\cite{6}. This was an
important achievement in the asymptotic theory.

In contrast to our paper, in A.~M.~Il'in's paper an arbitrary
smooth initial condition was considered. In~\cite{2}, in the case
$f(u)=u^2$, it was explained how the solution constructed there
can be used to obtain the global weak asymptotic for a more
general Cauchy problem. For this, it was proposed to consider an
interpolation of the initial function by linear splines. Also, we
can use different approach. For given smooth initial data
$u_0(x)$, $x\in {\bf R}$, we can find (assume finite) set of
points $x_k^0 \in {\bf R}$, $k\in {1,...,N}$, $N\in {\bf N}$,
which reaches the point of the gradient catastrophe in the moment
$t=t_k$. Now, instead of given initial data $u_0(x)$, $x\in {\bf
R}$, we impose initial data $u_{0\eps}(x)$, $x\in {\bf R}$, which
differs from the function $u_0(x)$, $x\in {\bf R}$, in the
intervals $[x_k^0-\eps^{\mu},x_k^0+\eps^{\mu}]$, $0<\mu<1$, $k\in
{1,...,N}$. In those intervals the function $u_{0\eps}(x)$, $x\in
{\bf R}$, has form (10). It is obvious that we have
$$||u_0(x)-u_{0\eps}(x)||_{L^1({\bf R})}={\mathcal O}(\eps^{\mu}).$$
Then, we solve the Cauchy problem corresponding to new initial
data $u_{0\eps}(x)$, $x \in {\bf R}$, using method of
characteristics in a way that for the "inserted" parts (the one in
the intervals $[x_k^0-\eps^{\mu},x_k^0+\eps^{\mu}]$, $0<\mu<1$,
$k\in {1,...,N}$) we use "the new characteristics" given by (11)
and for the rest of the function $u_{0\eps}(x)$~(~$\equiv u_0(x)$
for $x\notin [x_k^0-\eps^{\mu},x_k^0+\eps^{\mu}]$) we use ordinary
characteristics. This will be the subject of further
investigations.

\section{Description of the formula\\
for the weak asymptotic solution}

To construct a weak asymptotic solution
describing the passage from~(4) to~(5), we introduce some
auxiliary constructions.

We define a function $\xi(x_0)$ as a solution of the implicit
equation
\begin{equation}
U+u^0_0=u_1(x_0)+u_1(\xi(x_0)),
\end{equation}
which is solvable due to~(3).

Obviously,  $\xi:\,[a_2,a_1]\to[a_2,a_1]$
is a smooth isomorphism and $\xi(\xi(x_0))=x_0$.

We introduce the function $U_1(x_0,\rho)$, 
by setting
\begin{equation}
U_1(x_0,\rho)=B_2(\rho)u_1(x_0)+B_1(\rho)u_1(\xi(x_0)),
\end{equation}
where the functions $B_i(\rho)$, $i=1,2$, are defined in
Lemma~4.1 and the function $\rho=\rho(\tau)$ is defined below,
see~(12), (13) and 
$$
\tau=\frac{\varphi_{10}(t)-\varphi_{20}(t)}\ve,
\qquad
\varphi_{i0}(t)=a_i+f'(u_1(a_i))t, \ \ i=1,2.
$$

Note that, by (7), (8), and the formulas for $B_i$
at the end of Lemma~4.1,
we have
\begin{equation}
U+u^0_0-U_1(x_0,\rho)=U_1(\xi(x_0),\rho),
\quad
U_1(x_0,\rho)=u_1(x_0)+O(\rho^{-N}),\quad \rho\to\infty.
\end{equation}

We shall seek a weak asymptotic solution of problem (1)--(2) in
the form
\begin{align}
u_\ve(x,t)&=u^0_0
+\Big(U_1(x_0(x,t,\tau),\rho)-u^0_0\Big)\omega_1
\bigg(\frac{\phi_1-x}{\ve}\bigg)
\nonumber\\
&\qquad +\Big(U-U_1(x_0(x,t,\tau),\rho)\Big)\omega_2
\bigg(\frac{\phi_2-x}{\ve}\bigg),
\end{align}
where $\omega_i(z)\to0,1$ as $z\to\mp\infty$,
$\displaystyle\frac{d^\alpha\omega_i}{dz^\alpha}=O(|\tau|^{-N})$,
where $|z|\to\infty$, $\alpha>0$ and $N>0$ are arbitrary
numbers, and $\phi_i=\phi_i(t,\ve)$, $i=1,2$,
$x_0(x,t,\tau)$ are the desired functions.

As noted in the last Section (see Sec.~4.1),
the functions $\omega_i((\phi_i-x)/\ve)$
approximate (in the weak sense) the Heaviside functions
$H(\phi_i-x)$,
$$
\omega_i\bigg(\frac{\phi_i-x}\ve\bigg)
=H(\phi_i-x)+O_{\cD'}(\ve),\qquad i=1,2.
$$

We shall seek the functions $\phi_i=\phi_i(t,\ve)$, $i=1,2$, in the form
$$
\phi_i=\hat\varphi_i(t,\tau)+\psi_0\hat\phi(t,\tau),\qquad i=1,2,.
$$Here, $\hat\phi(t,\tau)$ is such that it satisfies
$\hat\phi(t,\tau)\big{|}_{\tau\to\infty}=0$. Furthermore,
$\hat\varphi_i(t,\tau)$ is an analog of the trajectories
of weak discontinuities of $\varphi_i(t)$
in~(4).
The functions $\hat\varphi_i(t,\tau)$, $i=1,2$, can be found from the
equations  for the "new characteristics"
and, as $\tau\to\infty$
(i.e., before the confluence of weak singularities),
these functions are close to the trajectories $\varphi_i(t)$
from the preceding section.

To find the functions $x_0(x,t,\tau)$,
we introduce the differential equation
for the "new characteristics"
\begin{equation}
\frac{dx}{dt}=B_2(\rho)f'(U_1(x_0,\rho))
+B_1(\rho)f'(U_1(\xi(x_0),\rho))
+q(\tau,\rho),\qquad
x\bigg|_{t=0}=x_0.
\end{equation}
The function $q(\tau,\rho)$ is assumed to be smooth and to
satisfy the estimate
\begin{equation}
|\tau q(\tau,\rho)|\leq\const.
\end{equation}

Its appearance itself is caused by the fact that the function
$U_1(x,\rho)$, which  replaces the function $u_1(x_0)$ in
formula~(4), depends on time
(via the function $\rho=\rho(\tau)$ determined in (13), (14)).
Therefore, this function is not preserved along the usual
trajectories corresponding to quasilinear equations.
The "new trajectories" are just given by Eq.~(11),
where the function~$\rho$ is determined as follows.
By $x(x_0,t,\tau)$ we denote the solution of~(11)
and introduce the functions $\hat\varphi_i(t,\tau)=x(a_i,t,\tau)$, $i=1,2$.
We set
$$
\rho=\frac{\hat\varphi_1(t,\tau)-\hat\varphi_2(t,\tau)}{\ve}
=\frac{\phi_1(t,\tau)-\phi_2(t,\tau)}{\ve}.
$$
We note that
$U_1(a_1,\tau)=B_2 u^0_0+B_1U$ and
$U_1(a_2,\tau)=B_2 U+B_1 u^0_0$;
hence from~(11) we easily obtain the following equation
for $\rho=\rho(\tau)$:
\begin{equation}
\frac{d\rho}{d\tau}=(B_2(\rho)-B_1(\rho))
(f'(B_2 u^0_0+B_1U)-f'(B_2 U+B_1 u^0_0))
(\psi'_0)^{-1}.
\end{equation}
Obviously, by definition,
\begin{equation}
\rho\tau^{-1}\to1\quad\text{as}\quad \tau\to\infty
\quad\text{and}\quad \frac{d\rho}{d\tau}>0\quad
(\text{since}\quad \psi'_0<0).
\end{equation}
We denote the right-hand side of (13) by $G(\rho)$.
Obviously, $G(\rho_0)=0$, where $\rho_0$ is a number
such that $B_1(\rho_0)=B_2(\rho_0)$, and hence (see Lemma~4.1)
\begin{equation}
B_1(\rho_0)=B_2(\rho_0)=1/2.
\end{equation}
We assume that~$\rho_0>0$. This is a condition imposed on~$B_j$.
It is easy to verify that
$\displaystyle\frac{dG}{d\rho}\bigg|_{\rho=\rho_0}=0$,
while
\begin{equation}
\frac{d^2G}{d\rho^2}\bigg|_{\rho=\rho_0}
=-8B^{\prime 2}_{2\rho}(U-u^0_0)f''\bigg(\frac{U+u^0_0}{2}\bigg)\ne0.
\end{equation}

It follows from Eq.~(13) and inequality~(16) that
the relations
\begin{equation}
\rho\to\rho_0+O(1/|\tau|),
\qquad
\dot\rho=O(1/|\tau|^2)
\end{equation}
hold as $\tau\to-\infty$.

Thus, independently of (11), the function $\rho=\rho(\tau)$ is
defined as a solution of problem (13), (14). Therefore, the
function $x(x_0,t,\tau)$ from~(11) is also defined. Now put
\begin{equation}
\hat{x}(x_0,t,\tau)=X(x_0,t)+\psi_0 X_1(x_0,\tau),
\end{equation}
where
$$
X(x_0,t)=x_0+f'(u_0(x_0))t=x_0\psi_0(\psi^0_0)^{-1}+bt.
$$
Inserting $\hat{x}$ in (11) instead of $x$ we have
\begin{align}
X_1
&=\frac1{\psi'_0\tau}\int^\tau_0\Big[
B_2(\rho)f'(U_1(x_0,\rho))+B_1(\rho)f'(U_1(\xi(x_0),\rho))
\nonumber\\
&\qquad
+q(\tau',\rho)-f'(u_0(x_0))\Big]\,d\tau'
\nonumber\\
&=\frac1{\psi'_0\tau}\int^\tau_0\Big[
B_2(\rho)f'(U_1(x_0,\rho))+B_1(\rho)f'(U_1(\xi(x_0),\rho))
\nonumber\\
&\qquad
+q(\tau',\rho)\Big]\,d\tau'
-(\psi^0_0)^{-1} x_0-b(\psi'_0)^{-1}.
\end{align}
It is easy to verify that the following representation is true:
$$
\hat{x}
=x^*+\frac{\psi_0}{\psi'_0\tau}\int^\tau_0
\Big[B_2(\rho)f'(U_1(x_0,\rho))+B_1(\rho)f'(U_1(\xi(x_0),\rho))
+q(\tau',\rho)\Big]\,d\tau',
$$
which follows from the identity
$$
x_0\psi_0(\psi^0_0)^{-1}+bt
-\psi_0(\psi^0_0)^{-1} x_0-\psi_0 b(\psi'_0)^{-1}
=-b\frac{\psi^0_0}{\psi'_0}=x^*.
$$
It is not difficult to see that the solution $\hat{x}$ given by
formula (19) is not the exact solution of (11). 
Actually, for $t=0$ (i.e.
for $\tau\to +\infty$)from (19) we obtain
$$(X_0+\psi_0X_1)|_{t=0}=x_0+{\mathcal O}(\eps).$$ 
Obviously, for
$t\in [0,T]$, $T\in {\bf R}$, we have:
$$
x(x_0,t,\tau)=\hat{x}(x_0,t,\tau)+O(\eps).
$$

It is easy to verify that the term $ O(\eps)$ 
in the last relation has the form
\begin{align*}
O(\eps)=\psi_0 X_1\Big|_{t=0}
&=\frac{\psi_0}{\psi'_0\tau}\int^\infty_0
\Big[B_2(\rho)f'(U_1(x_0,\rho))+B_1(\rho)f'(U_1(\xi(x_0),\rho))
\\
&\qquad 
+q(\tau',\rho)-f'(u_0(x_0))\Big]\,d\tau'
\od\ve g(x_0). 
\end{align*}
Finally, we obtain
$$
x(x_0,t,\tau)=\hat{x}(x_0,t,\tau)+\ve g(x_0).
$$

Let us calculate the derivative $\frac{\pa \hat{x}}{\pa x_0}$.
By (18), (19), we have
\begin{align}
\frac{\pa x}{\pa x_0}
&=\frac{\psi_0}{\psi'_0\tau}\int^\tau_0\frac{\pa }{\pa x_0}
\Big[B_2(\rho)f'(U_1(x_0,\rho))
+B_1(\rho)f'(U_1(\xi(x_0),\rho))\Big]\,d\tau
\nonumber\\
&=\frac{\psi_0}{\psi'_0\tau}\int^\tau_0
\Big[B_2(\rho)f''(U_1(x_0,\rho))
-B_1(\rho)f''(U_1(\xi(x_0),\rho))\Big] \frac{\pa U_1}{\pa
x_0}(x_0,\rho)\,d\tau
\end{align}
Here we used the relation
$$
\frac{\pa U_1}{\pa x_0}(x_0,\rho) =-\frac{\pa U_1}{\pa
x_0}(\xi(x_0),\rho),
$$
which follows from the definition of the function $U_1(x_0,\rho)$
in~(9).

 We agree that
the symbol~$\sim$ denotes the following equivalence relation
\begin{equation}
f\sim g\leftrightarrow \lim \frac fg=\const\ne0.
\end{equation}
Then, as $\tau\to-\infty$, we have
\begin{gather*}
\frac{\pa U_1}{\pa x_0}\sim B_1-\frac12\sim\frac1\tau,
\qquad
B_2-\frac12\sim\frac1\tau,
\\
U_1(x_0,\rho)\to\frac{U+u^0_0}2,\qquad
U_1(\xi(x_0),\rho)\to\frac{U+u^0_0}2.
\end{gather*}
Therefore,
$$
\Big[B_2(\rho)f''(U_1(x_0,\rho))
-B_1(\rho)f''(U_1(\xi(x_0),\rho))\Big]\sim\frac1\tau.
$$
Hence the integral in~(20) converges as $\tau\to-\infty$ and
$$
\frac{\pa \hat{x}}{\pa x_0}\sim\frac{\psi_0}{\psi'_0\tau},\qquad
\tau\to-\infty.
$$
As $\tau\to\infty$, we have $U_1(x_0,\rho)\to u_1(x_0)$
(since $B_2\to1$) and the integrand in~(20)
tends to the limit
$$
f''(u_1(x_0))\frac{\pa u_1}{\pa x_0}(x_0)=\frac{\psi'_0}{\psi^0_0}.
$$
Thus
$$
\frac{\pa \hat{x}}{\pa x_0}\to\frac{\psi_0}{\psi^0_0},\qquad
\tau\to\infty.
$$

We note that the solvability of the equation 
$x(x_0,t,x)=x$ with respect to~$x_0$ globally in~$t$
can hardly be ignored.

In our constructions, 
we shall hence use the following approximate expression for the
solution of Eq.~(11), namely,
$$
x(x_0,t,\tau)=\hat{x}(x_0,t,\tau)
+\ve (g(x_0)+Ax_0),
$$
where $A>0$, $A=\const$.
Clearly, we have
$$
x\Big|_{t=0}=x_0+\ve Ax_0,
$$
and hence $x_0(x,t,\tau)|_{t=0}=x-\ve Ax+O(\ve)$.

As is easy to verify, this means that, 
in the sense of $O_{\cal D'}$-estimates,  
the initial condition in (???) will be satisfied 
with accuracy up to $O_{\cal D'}(\ve)$.
We prove that the constant~$A$ can be chosen so that 
the inequality 
$$
\frac{\pa x}{\pa x_0}>0
$$
holds uniformly in~$t$.

We have
\begin{align*}
\frac{\pa x}{\pa x_0}
&=1-kt+\frac{\psi_0}{\psi'_0\tau}\int^\tau_0
\Big[B_2 f''(U_1(x_0,\rho))-B_1 f''(U_1(\xi(x_0),\rho))\\
&\qquad 
-f''(u_1(x_0))\frac{\pa u}{\pa x_0}\Big]\,d\tau'
+\ve g'(x_0)+\ve A.
\end{align*}
Recall that $t^*=k^{-1}$, $\psi_0(t^*)=0$, 
$\tau=\psi_0(t)/\ve$.
Hence for $t\leq t^*$, by Lemma~4.2, we have the estimate
$$
\frac{\psi_0}{\psi'_0\tau}\int^\tau_0
\Big[B_2 f''(U_1(x_0,\rho))-B_1 f''(U_1(\xi(x_0),\rho))
-f''(u_1(x_0))\frac{\pa u_1}{\pa x_0}\Big]\,d\tau'
=O(\ve)
$$
Similarly, for $t\geq t^*$, 
we have 
$$
1-kt+\frac{\psi_0}{\psi'_0\tau}\int^\infty_0
\Big[B_2 f''(U_1(x_0,\rho))-B_1 f''(U_1(\xi(x_0),\rho))
-f''(u_1(x_0))\frac{\pa u_1}{\pa x_0}\Big]\,d\tau'
=O(\ve).
$$
It follows from these estimates that there is a possibility 
to choose the constant~$A$. Thus the equation
\begin{equation}
X_0(x_0,t)+\psi_0 X_1(x_0,t,\ve)+\ve(g(x_0)+Ax_0)=x
\end{equation}
can be globally solved with respect to~$x_0$. 

In this case, the derivatives of the exact solution of Eq.~(11)
differ from the function in the right-hand side of (22)
and from the function $\hat{x}(x_0,t,\tau)$ 
by $O(\ve)$.
Therefore, in what follows, to simplify the calculations, 
we shall use all these functions.

\section{Construction\\ of the weak asymptotic solution}

We substitute the function $u_\ve(x,t)$ into Eq.~(1).
Using Lemma~4.1 and the formula for weak asymptotic
of the approximations in Sec.~4, we obtain
\begin{align}
\frac{\pa u_\ve}{\pa t}+\frac{\pa }{\pa x}f(u_\ve)
&=\phi_{1t}(U_1(x_0(\phi_1,t,\tau),\rho)-u^0_0)\delta(x-\phi_1)
\nonumber\\
&\quad
+\phi_{2t}(U-U_1(x_0(\phi_2,t,\tau),\rho))\delta(x-\phi_2)
\nonumber\\
&\quad
+\frac{\pa U_1}{\pa x_0}\frac{\pa x_0}{\pa t}
[H(\phi_1-x)-H(\phi_2-x)]
\nonumber\\
&\quad
+\frac{\pa U_1}{\pa x_0}
\frac{\pa x_0}{\pa x}
\big[B_2(\rho) f'(U_1(x_0(x,t,\tau),\rho))
\nonumber\\
&\quad
+B_1(\rho) f'(U_1(\xi(x_0(x,t,\tau)),\rho))\big]
[H(\phi_1-x)-H(\phi_2-x)]
\nonumber\\
&\quad
-\delta(x-\phi_1)
\Big[B_2(\rho)\Big( f(U_1(x_0(\phi_1,t,\tau),\rho))-f(u^0_0)\Big)
\nonumber\\
&\quad
+B_1(\rho)\Big( f(U)-f(U_1(\xi(x_0(\phi_1,t,\tau)),\rho))\Big)\Big]
\nonumber\\
&\quad
-\delta(x-\phi_2)
\Big[B_1(\rho)\Big( f(U_1(\xi(x_0(\phi_2,t,\tau)),\rho))-f(u^0_0)\Big)
\nonumber\\
&\quad
+B_2(\rho)\Big( f(U)-f(U_1(x_0(\phi_2,t,\tau),\rho))\Big)\Big]
\nonumber\\
&\quad
+\frac{\pa U_1}{\pa t}(x_0,\rho)\bigg|_{x_0=x_0(x,t,\tau)}
[H(\phi_1-x)-H(\phi_2-x)]+O_{\cD'}(\ve).
\end{align}

Although there are rather many terms
on the right-hand side,
it is easy to understand this formula.
The terms containing the factors $(H(\phi_1-x)-H(\phi_2-x))$
correspond to the substitution into the equation
of the function $u_\ve$ determined in~(10)
between the points $x=\phi_i$
and with Lemma~4.1 taken into account.

The terms containing the delta functions, i.e., the factors
$\delta(x-\phi_i)$, $i=1,2$, appear due to the fact that
$U(x_0(\phi_i(t,\tau),\rho))\ne u^0_0$  and
$U(x_0(\phi_i(t,\tau),\rho))\ne U$,
but as $\tau\to\infty$ (i.e., before the interaction)
we have $\rho\sim\tau$ (see (14)) and hence, for any $N>0$,
we have
$$
U(x_0(\phi_1(t,\tau),\rho))-u^0_0=O(\ve^N),\qquad
U(x_0(\phi_2(t,\tau),\rho))-U=O(\ve^N).
$$
We start analyzing the terms in (22) from the last one
(which has the estimate $O(\ve^{-1})$ in the $C$-norm):
\begin{align}
\frac{\pa U_1}{\pa t}
&=\ve^{-1}\psi'_0\dot\rho
\big[B'_{2\rho}u_1(x_0)+B'_{1\rho}u_1(\xi(x_0))\big]
\nonumber\\
&=\ve^{-1}\psi'_0\dot\rho
B'_{2\rho}\big[U+u^0_0-2u_1(x_0(x,t,\tau))\big].
\end{align}
Applying the Taylor formula at the points $x=\phi_1$ and
$x=\phi_2$, for any test function $\eta(x)$,
we obtain
\begin{align}
&\ve^{-1}\psi'_0\dot\rho B'_{2}
\int^{\phi_1}_{\phi_2}[U+u^0_0]\eta(x)\,dx
\nonumber\\
&\qquad
=\frac{\psi'_0}2\rho \dot\rho B'_{2\rho}(U+u^0_0)
\auf\{\delta(x-\phi_1)+\delta(x-\phi_2)\},\eta(x) \zu+O(\ve),
\end{align}
and $B'_{2\rho}=O(|\rho|^{-N})$ for any $N>0$ as $\rho\to\infty$,
$B'_{2\rho}\to\const$ as $\rho\to\rho_0$ ($\tau\to-\infty$),
and $\dot\rho=O(|\tau|^{-2})$ as $\tau\to-\infty$ (see (17)).

Let us consider the remaining term.
We have
\begin{equation}
\int^{\phi_1}_{\phi_2} u_1(x_0(x,t,\tau))\eta(x)\,dx
=\int^{a_1}_{a_2}u_1(x_0)\eta(x(x_0,t,\tau)+\psi_0\hat\phi)
\frac{dx}{dx_0}\,dx_0.
\end{equation}
Let us note that (see (18))
\begin{equation}
\frac{\pa x}{\pa x_0}=\psi_0
\bigg((\psi^0_0)^{-1}+\frac{\pa X_1}{\pa x_0}\bigg)
=\ve\tau
\bigg((\psi^0_0)^{-1}+\frac{\pa X_1}{\pa x_0}\bigg).
\end{equation}
Hence the right-hand side in (23) is bounded in the weak sense
as $\ve\to0$.
We now note that the following relations hold:
\begin{gather*}
\eta(x(x_0,t,\tau)+\psi_0\hat\phi)
=\eta(x(a_i,t,\tau)+\psi_0\hat\phi)
+\eta'_x\frac{\pa x}{\pa x_0}\bigg|_{x_0=c_i},\\
c_i\in(a_i,x_0),\qquad i=1,2.
\end{gather*}
Recalling that $x(a_i,t,\tau)=\hat{\varphi}_i(t,\tau)$, $i=1,2$,
and again using (26), (14) and (17) we obtain
\begin{align*}
&\ve^{-1}\psi'_0\dot\rho B'_{2}
\int^{\phi_1}_{\phi_2} u_1(x_0(x,t,\tau))\eta(x)\,dx
\\
&\qquad
=\frac12\auf\delta(x-\phi_1)+\delta(x-\phi_2),\eta(x)\zu
\psi'_0\dot\rho\tau B'_{2}\\
&\qquad\qquad
\times
\int^{a_1}_{a_2} u_1(x_0)
\bigg((\psi^0_0)^{-1}+\frac{\pa X_1}{\pa x_0}\bigg)\,dx_0
+O(\ve).
\end{align*}

Finally, we have
$$
\frac{\pa U_1}{\pa t}[H(\phi_1-x)-H(\phi_2-x)]
=g(\tau,\rho)\big(\delta(x-\phi_1)+\delta(x-\phi_2)\big)
+O_{\cD'}(\ve),
$$
where
\begin{equation}
g(\tau,\rho)=\frac{\psi'_0}2\rho\dot\rho B'_{2\rho}(U+u^0_0)
-\psi'_0\dot\rho\tau B'_{2\rho}
\int^{a_1}_{a_2} u_1(x_0)
\bigg((\psi^0_0)^{-1}+\frac{\pa X_1}{\pa x_0}\bigg)\,dx_0.
\end{equation}

It is easy to see that, by formula~(17),
we have the estimate
$$
|\tau^2 g(\tau,\rho)|\leq\const.
$$
Moreover, the function $g(\tau,\rho)$ is integrable,
and the integral $\int^\tau_0 g(\tau,\rho)\,d\tau$ converges.
Indeed,
the integral of the first term converges because of the estimates
given after formula~(24),
and the integral of the second term, in its properties,
coincides with the last integral in formula~(20).

Now we consider the remaining terms that contain the difference
$H(\phi_1-x)-H(\phi_2-x)$ as the multiplier.
For any function $\eta(x)\in C^\infty_0$,
taking into account the relation
$$
\frac{\pa x_0}{\pa t}=-\frac{\pa x_0}{\pa x}\frac{\pa x}{\pa t},
$$
we have
\begin{align}
&\bigg\langle\bigg[
\frac{\pa U_1}{\pa x_0}\frac{\pa x_0}{\pa t}
+\frac{\pa U_1}{\pa x_0}\frac{\pa x_0}{\pa x}
\Big(B_2(\rho) f'(U_1(x_0(x,t,\tau),\rho))
\nonumber\\
&\qquad
+B'_1(\rho) f'(U+u^0_0-U_1(x_0(x,t,\tau),\rho))\Big)\bigg]
(H(\phi_1-x)-H(\phi_2-x)),\eta(x)\bigg\rangle
\nonumber\\
&\quad
=\int^{\phi_1}_{\phi_2}
\frac{\pa U_1}{\pa x_0}\frac{\pa x_0}{\pa x}
\bigg[-\frac{\pa x}{\pa t}+B_2(\rho)f'(U_1(x_0(x,t,\tau),\rho))
\nonumber\\
&\qquad
+B_1(\rho)f'(U_1(\xi(x_0(x,t,\tau)),\rho))\bigg]\eta(x)\,dx.
\end{align}

By (11), the expression in square brackets
on the right-hand side of (28) is just $q(\tau,\rho)$.

We consider the integral
$$
\int^{\phi_1}_{\phi_2}\frac{\pa U_1}{\pa x_0}
\frac{\pa x_0}{\pa x}\eta(x)\,dx
$$
and pass to the variables~$x_0$ precisely as in (25).
We obtain
$$
\int^{\phi_1}_{\phi_2}\frac{\pa U_1}{\pa x_0}
\frac{\pa x_0}{\pa x}\eta(x)\,dx
=\int^{a_1}_{a_2}\frac{\pa U_1}{\pa x_0}
\eta(x(x_0,t,\tau)+\psi_0\hat\phi)\,dx_0.
$$
Recall that
$$
\frac{\pa U_1}{\pa x_0}\sim\frac1\tau,
\quad \tau\to-\infty,\qquad
\frac{\pa U_1}{\pa x_0}\to\frac{\pa u_1(x_0)}{\pa x_0},
\quad \tau\to\infty.
$$
From the conjectural estimate~(12) for the function
$q(\tau,\rho)$,
using the Taylor formula as in (24), we obtain
$$
q\int^{\phi_1}_{\phi_2}\frac{\pa U_1}{\pa x_0}
\frac{\pa x_0}{\pa x}\eta(x)\,dx
=\frac q2\int^{a_1}_{a_2}\frac{\pa U_1}{\pa x_0}\,dx_0
\auf\delta(x-\phi_1)+\delta(x-\phi_2),\eta(x)\zu
+O(\ve).
$$

Taking into account the definition of the function
$U_1(x_0,\rho)$,
we can easily calculate
the integral on the right-hand side of the last formula
and obtain
\begin{align}
&q\int^{\phi_1}_{\phi_2}\frac{\pa U_1}{\pa x_0}
\frac{\pa x_0}{\pa x}\eta(x)\,dx
\nonumber\\
&\qquad
=\frac {q(B_2-B_1)(u^0_0-U)}2
\auf\delta(x-\phi_1)+\delta(x-\phi_2),\eta(x)\zu
+O(\ve).
\end{align}

We choose the function $q(\tau,\rho)$
so that the following relation hold:
\begin{equation}
q(B_2-B_1)(u^0_0-U)=-g(\tau,\rho).
\end{equation}

Obviously, we have
\begin{gather*}
q(\tau,\rho)\sim g(\tau,\rho)=O(\tau^{-N})\quad \forall N,
\quad\tau\to\infty;
\\
q(\tau,\rho)\sim \tau g(\tau,\rho)\sim\frac1\tau,\quad\tau\to-\infty.
\end{gather*}
Hence the estimate~(12) holds and our constructions
that lead to (29) are well defined.

It is left to obtain the function $\hat{\phi}$ appearing in the definition of the 
functions $phi_i$, $i=1,2$. To do that we will use the results from Section 4.3.
Equating with zero the remaining coefficients of
$\delta(x-\phi_i)$, $i=1,2$,
(only such expressions $\pmod{O_{\cD'}(\ve)}$ are left
on the right-hand side of~(22)),
we obtain
\begin{align}
&\phi_{1t}( U( x_0(\phi_1,t,\tau),\rho)-u^0_0)
-B_2(\rho)\big( f(U_1(x_0(\phi_1,t,\tau),\rho))-f(u^0_0)\big)
\nonumber\\
&\qquad
-B_1(\rho)\big( f(U)-f(U_1(\xi (x_0(\phi_1,t,\tau)),\rho)) \big)=0,
\\
&\phi_{2t} (U- U_1( x_0(\phi_2,t,\tau),\rho))
-B_2(\rho)\big( f(U)-f(U_1(x_0(\phi_2,t,\tau),\rho))\big)
\nonumber\\
&\qquad
-B_1(\rho)\big(f (U_1(\xi(x_0(\phi_2,t,\tau)),\rho))-f(u^0_0)\big)=0.
\end{align}According to (47) we have to prove that preceding equations
are correct when $\tau \to \pm\infty$ and to find $\hat{\varphi}$ such that 
their sum be equal to zero.

By the definition of the functions $\phi_i(t,\tau)$, $i=1,2$,
as $\tau\to\infty$ (i.e., before the interaction),
the limit of the expressions on the left-hand side
of relations (29), (30) is equal to zero,
and these relations admit the estimate $O(\tau^{-N})$
for any $N>0$ as $\tau\to\infty$.
This follows from the relations:
$\rho/\tau\to1$ as $\tau\to\infty$ and
\begin{equation}
B_2=1+O(\rho^{-N}),\qquad B_1=O(\rho^{-N}),\qquad
\rho\to\infty\qquad (\tau\to\infty).
\end{equation}

We write the limit of these relations for $\tau\to-\infty$.
Recall that
\begin{equation}
B_i(\rho)=\frac12+O(|\tau|^{-1}),\qquad
\rho=\rho_0+O(|\tau|^{-1}),\qquad
\tau\to-\infty,\qquad i=1,2.
\end{equation}

Therefore, denoting the limit of $\phi_{it}$ as $\tau\to-\infty$
by $\phi^-_{it}$, $i=1,2$, we obtain
\begin{equation}
\phi^-_{it}\bigg(\frac{U-u^0_0}2\bigg) =\frac12 (f(U)-f(u^0_0)),
\qquad i=1,2,
\end{equation}
or
\begin{equation}
\phi^-_{1t}=\frac{f(U)-f(u^0_0)}{U-u^0_0} =\phi^-_{2t}.
\end{equation}

Denoting, as usual,
$$
\frac{f(U)-f(u^0_0)}{U-u^0_0}=\frac{[f]}{[u]},
$$
we can determine
the general limit $\phi^-(t)$ of the functions $\phi_i(\tau,t)$, $i=1,2$,
as $\tau\to-\infty$ by the relation
\begin{equation}
\phi^-=\phi^-(t^*)+\frac{[f]}{[u]}(t-t^*).
\end{equation}

Relations (36) (or (38)) mean that, for $t>t^*$,
the trajectories $x=\phi_{1}$ and $x=\phi_2$ are close
to the line
$$
x-x^*=\frac{f(U)-f(u^0_0)}{U-u^0_0} (t-t^*),
$$
i.e.,
to the trajectory of the shock wave (5).

Let us investigate the trajectories $x=\phi_i$, $i=1,2$,
in more detail.

By $\omega(z)$ we denote the function satisfying
the same conditions as the functions $\omega_i$, $i=1,2$, in (10).

We prove that the following relations hold:
\begin{equation}
\phi_i(t,\tau)-\check\phi_i(t,\tau)=O(\ve),\qquad i=1,2,
\end{equation}
where
$$
\check\phi_i(t,\tau)=(1-\omega(\tau))\hat\varphi_i(\tau,t)
+\omega(\tau)\bigg(x^* + \frac{[f]}{[u]}(t-t^*)\bigg).
$$
Here $\hat\varphi_i(\tau,t)=X(a_i,t)+\psi_0 X_1(a_i,\tau)$,
(see (18)),
$x^*=\varphi_{10}(t^*)=\varphi_{20}(t^*)$,
and $\phi_i(t,\tau)$, $i=1,2$, are the desired trajectories of
singularities determined $\pmod{O(\ve)}$ by the relations
$$
\phi_i(t,\tau)=\hat\varphi_i(\tau,t)+\psi_0 \hat\phi.
$$

To prove (39), it suffices to set $\phi^-(t^*)=x^*$ in (38) 
and
to note that $\hat\varphi_i(0,t^*)=x^*$, $i=1,2$.
It remains to note that the functions $\hat\varphi_i(\tau,t)$
can be represented in the form
\begin{equation}
\hat\varphi_i(\tau,t)
=x^*+\psi_0(\psi'_0\tau)^{-1}\int^\tau_0
[B_2 f'(U(a_i,\rho))+B_1 f'(U(\xi(a_i),\rho))+q(\tau',\rho)]\,d\tau'.
\end{equation}
This follows from (18), (19) with the relation 
$x^*=b t^*=-b\psi^0_0/\psi'_0$ taken into
account (see formula~(3)).

Now we apply Lemma~4.2 and see that relation~(39) is proved.
The statement we have proved means that $\check\phi_i$
from (39)
provide a family of expressions for trajectories close to those
trajectories we want to construct.
These approximate trajectories, with accuracy~$O(\ve)$,
are independent of the choice of the function $\omega(\tau)$.
It is only required that this function satisfy
same conditions as the functions $\omega_i$, $i=1,2$, from~(10).

Now let us calculate the function $x_0(\phi_i,t,\tau)$, $i=1,2$.
By definition, this is the initial point of the trajectory
$x=\phi_i(t,\tau)$, $i=1,2$.
Clearly, for $t<t^*$,
we have
$\phi_i(t,\tau)=\hat\varphi_i(t,\tau)+O(\ve)$ and
$x_0(\varphi_i,t,\tau)=a_i$.
For $t>t^*$, we have $\phi_i(t,\tau)-\phi^-(t)\to0$ as
$\ve\to0$.
By relation (37), for $\phi^-(t^*)=x^*$,
we see that in this case  the initial point is
$$
\check x=\phi^-(0)=x^*-\frac{[f]}{[u]}t^*.
$$
By the inequalities $f'(U)<[f]/[u]<f'(u^0_0)$,
this implies that $\check x\in(a_2,a_1)$.

We set
$$
\hat{X}_0(\phi_i,\tau)=a_i+\Omega(\tau)(\check x-a_i),\qquad i=1,2,
$$
where $\Omega(\tau)$ is some (generating) function satisfying same
conditions as the functions $\omega_i$, $i=1,2$, from~(10).

Let us prove the relations
\begin{equation}
\phi_i(t,\tau)-\big(x(\hat{X}_0(\phi_i,\tau),t,\tau)
+\psi_0\hat\phi\big)=O(\ve),
\qquad i=1,2.
\end{equation}
We restrict ourselves only to the case $i=1$.
We have
\begin{align*}
U_1(a_1,\rho)
&=U_1(\hat{X}_0,\rho)-\Omega(\check x-a_i)\frac{\pa U_1}{\pa x_0}
\big(a_1+\alpha\Omega(\check x-a_i),\rho\big),
\\
U_1(\xi(a_1),\rho)
&=U+u^0_0-U_1(a_1,\rho)\\
&=U+u^0_0-U_1(\hat{X}_0,\rho)-(\check x-a_i)\frac{\pa U_1}{\pa x_0}
\big(a_1+\alpha\Omega(\check x-a_i),\rho\big),
\end{align*}
where $\alpha\in(0,1)$.

From these relations, formula~(40), 
and representation for $\hat{x}$ from Section~2, 
we obtain
\begin{align}
&\hat\varphi_1(t,\tau)-\hat{x}(\hat{X}_0(\phi_1,\tau),t,\tau)
\\
&\qquad =(\check x-a_i)\psi_0(\psi'_0\tau)^{-1}
\int^\tau_0
\bigg[B_2\frac{\delta f'}{\delta u}
(U_1(a_1,\rho);U_1(\hat{X}_0,\rho))
\nonumber\\
&\qquad\qquad
-B_1\frac{\delta f'}{\delta u}
(U+u^0_0-U_1(a_1,\rho);U+u^0_0-U_1(\hat{X}_0,\rho))\bigg]
\nonumber\\
&\qquad\qquad \times
\Omega\frac{\pa U_1}{\pa x_0}
\big(a_1+\alpha\Omega(\check x-a_i),\rho\big)\,d\tau'
\nonumber
\end{align}
where
$$
\frac{\delta f'}{\delta u}(A,B)
=\frac{f'(A)-f'(B)}{A-B}\underset{A\to B}{\longrightarrow} f''(A).
$$
We note that the integral on the right-hand side of (42)
converges as $\tau\to+\infty$
because the function $\Omega$ is contained
in the integrand.
The convergence of the integral as $\tau\to-\infty$
can be verified in the same way as the convergence
of the last integral on the right-hand side of (20).
Hence, by Lemma~4.2, we have
$$
\hat\varphi_1(t,\tau)-x\big(\hat{X}_0(\phi_1,\tau),t,\tau\big)
=O(\ve),
$$
and hence, by (38) and (39), we obtain~(41).
From~(41) we obtain the relation
\begin{equation}
U_1(x_0(\phi_i,t,\tau),\rho)-U_1(\hat{X}_0(\phi_i,\tau),\rho)
=O(\ve),
\qquad i=1,2.
\end{equation}

By construction, the limits of the expressions on the left-hand
sides in (31) and (32) are equal to zero
as $\tau\to\infty$ (i.e., before the interaction).
Moreover,
the difference between the limit and the prelimit
expression is $O(\rho^{-N})=O(\tau^{-N})$ for any $N>0$.

By (32), these expressions also tend to zero
as $\tau\to-\infty$, and the difference between the limit
and the prelimit expression is
$O(B_1-1/2)=O(\rho-\rho_0)=O(|\tau|^{-1})$, $\tau\to-\infty$.
Therefore,
by the results of Sec.~4.2 about the linear independence,
for the sum of terms with $\delta$-functions in (22)
to admit the estimate $O_{\cD'}(\ve)$,
it is sufficient that the sum of expressions
on the left-hand sides of (31) and (32) be equal to zero.
Thus we obtain the equation
\begin{align}
&\phi_{2t}(U-U_{1(2)})
+\phi_{1t}(U_{1(1)}-u^0_0)
\nonumber\\
&\qquad
= B_2(\rho)\big(f(U_{1(1)})-f(u^0_0)\big)
+B_1(\rho)\big(f(U)-f(\hat U_{1(1)})\big)
\nonumber\\
&\qquad\qquad
+B_2(\rho)\big(f(U)-f(U_{1(2)})\big)
+B_1(\rho)\big(f(\hat U_{1(2)})-f(u^0_0)\big).
\end{align}
Here, for brevity, we denote
$$
U_{1(i)}=U_1(x_0(\phi_i,t,\tau),\rho),\qquad
\hat U_{1(i)}=U_1(\xi(x_0(\phi_i,t,\tau)),\rho), \ \ i=1,2.
$$
We note that
$$
\phi_{it}=\hat\varphi_{it}+\psi'_0\frac{d}{d\tau}(\tau\hat\phi),
\qquad i=1,2,
$$

We agree to denote $f\approx g$ if
$$
\lim \frac fg=1.
$$
It is easy to verify that as $\tau\to\infty$, we have
\begin{equation}
U-U_{1(2)}
\approx U_{1(1)}-u^0_0
\approx U-\hat{U}_{1(1)}
\approx \hat U_{1(2)}-u^0_0
\approx B_1(U-u^0_0).
\end{equation}
Similarly,
\begin{align}
f(U_{1(1)})-f(u^0_0)\approx f'(u^0_0) B_1(U-u^0_0),
\\
f(U)-f(U_{1(2)})\approx f'(U) B_1(U-u^0_0),
\nonumber\\
f(U)-f(U_{1(1)})\approx f'(U) B_1(U-u^0_0),
\nonumber\\
f(\hat U_{1(2)})-f(u^0_0)\approx f'(u^0_0) B_1(U-u^0_0).
\nonumber
\end{align}

Next, by~(16), we have $B'_2(\rho)\sim1-B_2$
and hence the relation $g\sim 1-B_2$ holds as $\tau\to+\infty$.

As $\tau\to-\infty$,
the coefficient of $\frac{d}{d\tau}(\tau\hat{\phi})$
in Eq.~(44) is equal to $U-u^0_0\ne0$.
Therefore, Eq.~(44) is solvable for $\hat \phi$
and its solution is a bounded function decreasing as
$\tau\to\infty$.

To write the solution of Eq.~(44), we note that,
with accuracy $O(\ve)$,
by (41), we can replace the arguments $x_0(\phi_i,t,\tau)$
by $X_0(\phi_i,\tau)$ in the functions $U_{i(j)}$,
and by (38), the function $X_0(\phi_i,\tau)$
can be determined actually independent of the
functions~$\phi_i$ (everywhere here $i,j=1,2$).
Hence Eq.~(44) is indeed a linear equation with respect
to~$\hat\phi$ and its solution can be easily found.

This solution has the form
\begin{align*}
\hat\phi&=(\psi'_0\tau)^{-1}\int^\tau_0
\big(U-u^0_0-U_1(\hat{X}_0(\phi_2,\tau),\rho)
+U_1(\hat{X}_0(\phi_1,\tau),\rho)\big)^{-1}\\
&\qquad
\times
\Big(-\hat\varphi_2[U-U_1(\hat{X}_0(\phi_2,\tau),\rho)]
-\hat\varphi_1[U_1(\hat{X}_0(\phi_1,\tau),\rho)-u^0_0]\\
&\qquad
+\big\{B_2(\rho)(f(U_{1(1)})-f(u^0_0))
+B_1(\rho)(f(U)-f(\hat U_{1(1)}))\\
&\qquad
+B_2(\rho)(f(U)-f(U_{1(2)}))
+B_1(\rho)(f(\hat U_{1(2)})-f(u^0_0))\big\}\Big)\,d\tau'.
\end{align*}
By (45) and (46), the integral on the right-hand side in the
last relation converges as $\tau\to\infty$ and
$\hat\phi=O(\tau^{-1})$ as $\tau\to\infty$.

\section {Auxiliary formulas and statements of weak asymptotic
method}

\subsection{Nonlinear superposition of approximations of
Heaviside functions}

Suppose that $\omega_j(z)\to0,1$ as $z\to-\infty$, $z\to\infty$,
$\displaystyle \frac{d^\alpha\omega_j}{dz^\alpha}=O(|z|^{-N})$, $j=1,2$,
$|z|\to\infty$, $N$ is a sufficiently large number,
and $\varphi_1$, $\varphi_2$ are some continuous functions of
the variable~$t$.

It is easy to verify that the functions
$\omega_j((x-\phi(t))/\ve)$ approximate in the weak sense
the Heaviside function $H(x-\phi(t))$.
Indeed, the properties of the functions $\omega_j(z)$
imply the relations
$$
\omega_j(z)-H(z)=O(|z|^{-N}),\qquad N>0,\quad j=1,2.
$$
Hence, for any test function $\psi(x)$, we have
$$
\Big\langle \omega_j\Big(\frac{x-\phi}{\ve}\Big)
-H(x-\varphi),\psi\Big\rangle
=\ve\int \big(\omega_j(z)-H(z)\big)\psi(\phi+\ve z)\,dz=O(\ve).
$$

\begin{lemma}
For any $C^1$-function $f(x)$, the following relation
holds{\rm:}
\begin{align*}
&f\bigg(a+b\omega_1\bigg(\frac{\varphi_1-x}{\ve}\bigg)
+b\omega_2\bigg(\frac{\varphi_2-x}{\ve}\bigg)\bigg)
\\
&\quad
=f(a)
+H(\varphi_1-x)\{B_2(f(a+b)-f(a))+B_1(f(a+b+c)-f(a+c))\}
\\
&\quad
+H(\varphi_2-x)\{B_1(f(a+c)-f(a))+B_2(f(a+b+c)-f(a+b))\}
+O_{\cD'}(\ve),
\end{align*}
where
\begin{gather*}
B_j=B_j\bigg(\frac{\varphi_1-\varphi_2}{\ve}\bigg), \; j=1,2 \;\qquad
B_1+B_2=1,\\
B_2(z)\to1\quad\text{as}\quad z\to\infty,
\qquad
B_2(z)\to0\quad\text{as}\quad z\to-\infty.
\end{gather*}
\end{lemma}

\proof
First, we prove the relation
$$
f\bigg(a+b\omega_1\bigg(\frac{\varphi_1-x}{\ve}\bigg)
+c\omega_2\bigg(\frac{\varphi_2-x}{\ve}\bigg)\bigg)
=
f(a+bH(\varphi_1-x)+cH(\varphi_2-x))+O_{\cD'}(\ve).
$$
Indeed, we have
\begin{align*}
&f(a+b\omega_1+c\omega_2)
=f\Big(a+bH(\varphi_1-x)+c(\varphi_2-x)\Big)
\\
&\qquad
+f'\Big(a+\xi(b(\omega_1-H(\varphi_1-x)))+c(\omega_2-H(\varphi_2-x))\Big)
\\
&\qquad\qquad
\times
[(\omega_1-H(\varphi_1-x))b+(\omega_2-H(\varphi_2-x))c].
\end{align*}

Now we verify that if $g(x,\varphi,\ve)$ is a bounded function,
then
$$
g(x,\varphi,\ve)[\omega_1-H(\varphi_1-x)]=O_{\cD'}(\ve).
$$
For any test function $\psi(x)$, we have
\begin{align*}
&\bigg|\int g(x,\varphi,\ve)
\bigg[\omega\bigg(\frac{\varphi_1-x}{\ve}\bigg)-H(\varphi_1-x)\bigg]
\psi(x)\,dx\bigg|\\
&\qquad
=\bigg|\ve\int g(\varphi_1+\ve z,\varphi,\ve)
[\omega(z)-H(z)]\psi(\varphi-\ve z)\,dz\bigg|\\
&\qquad
\leq \ve\const\int|\omega(z)-H(z)|\,dz.
\end{align*}
This implies
$$
f(a+b\omega_1+c\omega_2)
=f(a+bH(\varphi_1-x)+cH(\varphi_2-x))+O_{\cD'}(\ve).
$$
Next, it is easy to verify the relation
\begin{align*}
f(a+bH_1+cH_2)
&=f(a)+H_1[f(a+b)-f(a)] +H_2[f(a+c)-f(a)]\\
&\qquad
+H_1H_2\big(f(a+b+c)-f(a+c)-f(a+b)+f(a)\big),\\
&
H_j\overset{\rm def}{=} H(\varphi_j-x),\quad j=1,2.
\end{align*}
It remains to note that we have
\begin{gather*}
H(\varphi_1-x)H(\varphi_2-x)
=B_1H(\varphi_1-x)+B_2H(\varphi_2-x)
+O_{\cD'}(\ve),\\
B_1=\int \dot\omega_1(z)
\omega_2\bigg(z-\frac{\varphi_1-\varphi_2}{\ve}\bigg)\,dz,
\qquad B_2=1-B_1.
\end{gather*}
For the proof of these and similar relations,
see \cite{1,2,7}.
The proof of the lemma is complete.
\endproof

\subsection{Asymptotic linear independence}
If we want to consider linear combinations of generalized
functions with accuracy $O_{\cD'}(\ve^\alpha)$,
then we need to modify the notion of linear independence.
This modification plays the key role in considerations
related to the soliton interaction problem.

Indeed, let $\phi_1\ne \phi_2$ be independent of~$x$.
We consider the relation
$$
g_1\delta(x-\phi_1)+g_2\delta(x-\phi_2)
=O_{\cD'}(\ve^\alpha),\qquad \alpha>0,
$$
where $g_i$ are independent of $\ve$.
Obviously, we obtain the relations
$$
g_i=O_{\cD'}(\ve^\alpha),\qquad i=1,2,
$$
which, by our assumption, imply
$$
g_i=0,\qquad i=1,2.
$$

Everything is different~\cite{1}
if we assume that the coefficients~$g_i$, $i=1,2$, can depend on~$\ve$.
Here we consider only a special case of such dependence, which
we shall use later.
Namely, let
$$
g_i=A_i+S_i({\Delta \phi}/{\ve}), \qquad i=1,2,
$$
where $A_i$ are independent of $\ve$ and $S_i(\rho)$ decrease
as $|\rho|\to\infty$.

We assume that the estimate holds:
$$
|\rho S_i(\rho)|\leq\const,\qquad i=1,2.
$$

Let us find out what properties of the coefficients $g_i$
follow from the relation
$$
g_1\delta(x-\phi_1)+g_2\delta(x-\phi_2)=O_{\cD'}(\ve).
$$
Applying both sides of the equality
to a test function~$\varphi$,
we obtain
$$
g_1\varphi(\phi_1)+g_2\varphi(\phi_2)=O(\ve)
$$
or, which is the same,
\begin{equation}
[A_1\varphi(\phi_1)+A_2\varphi(\phi_2)]
+[S_1\varphi(\phi_1)+S_2\varphi(\phi_2)]
=O(\ve).
\end{equation}

Let us consider the expression in the second brackets.
Using Taylor's formula, we obtain
$$
[S_1\varphi(\phi_1)+S_2\varphi(\phi_2)]
=
S_1\varphi(\phi_1)+S_2\varphi(\phi_1)
+S_2(\phi_2-\phi_1)\varphi'(\phi_1+\theta \phi_2), \quad 0<\theta<1.
$$
Now we see that
$$
S_2({\Delta \phi}/{\ve})(\phi_2-\phi_1)
=\{-\rho S_2(\rho)\}\big|_{\rho=\Delta \phi/\ve}\cdot \ve=O(\ve),
$$
since the function $\rho S_2(\rho)$ is bounded uniformly
in $\rho\in\mathbb{R}^1$.

So we can rewrite relation (37) as
$$
A_1\varphi(\phi_1)+A_2\varphi(\phi_2)+(S_1+S_2)\varphi(\phi_1)=O(\ve).
$$
Hence, as the coefficients $A_i$ are independent of $\ve$,
we, as usual, obtain
\begin{equation}
A_1=0,\qquad A_2=0,\qquad S_1+S_2=0.
\end{equation}
Another method for analyzing relation (37) is the following.
We assume that $\phi_i(t)$ are smooth functions,
the relation $\phi_1(t^*)=\phi_2(t^*)$ holds for some $t=t^*$,
and, moreover, $\phi'_1(t^*)\not=\phi'_2(t^*)$.
Then
\begin{align*}
&\langle S_1\delta(x-\phi_1),\varphi\rangle
+\langle S_2\delta(x-\phi_2),\varphi\rangle
=S_1\varphi(x^*)+S_2\varphi(x^*)\\
&\qquad
+S_1 O(t-t^*)+S_2 O(t-t^*),\qquad
x^*=\phi_1(t^*)=\phi_2(t^*).
\end{align*}
But $O(t-t^*)\sim O(\Delta\phi)$. Therefore, we have
\begin{align*}
&(A_1+S_1)\delta(x-\phi_1)+(A_2+S_2)\delta(x-\phi_2)\\
&\qquad
=
A_1\delta(x-\phi_1)+A_2\delta(x-\phi_1)
+(S_1+S_2)\delta(x-x^*)+O_{\cD'}(\ve).
\end{align*}
We again obtain relations~(47).

\subsection{Complex germ lemma}
In this section,
in the form convenient for us,
we present the statement
that plays an important role
in Maslov's complex germ theory~\cite{8,9}.

\begin{lemma}
Let $f(t)\in C^1$, $f(t_0)=0$, and $f'(t_0)\ne0$.
Let $g(\tau,t)$ be a function that locally uniformly satisfies
the estimates
$$
|\tau g(\tau,t)|\leq\const,\qquad
|\tau g'_t(\tau,t)|\leq\const,\quad -\infty<\tau<\infty,
$$
and $g(\tau,t_0)=0$.
Then the inequality
$$
\bigg|g\bigg(\frac{f(t)}{\ve},t\bigg)\bigg|\leq C_T\ve,
$$
where $C_T=\const$,
holds on any interval $0\leq t\leq T$ that does not contain
zeros of the function $f(t)$ except~$t_0$.
\end{lemma}

\proof
The fraction $f(t)/(t-t_0)$ is locally bounded in~$t$.
The fraction $\tau g(\tau,t)/$ $(t-t_0)$ is also locally bounded.
We have
$$
g\bigg(\frac{f(t)}{\ve},t\bigg)
=\ve\bigg[g\bigg(\frac{f(t)}{\ve},t\bigg)(t-t_0)^{-1}\bigg]
\frac{f(t)}{\ve}\cdot\frac{t-t_0}{f(t)}.
$$
By the assumptions of the lemma,
on the interval under study,
the last multiplier on the right-hand side is bounded,
while the product of the intermediate multiplier
and the expression in square brackets is bounded
in view of the properties of the function $g(\tau,t)$.
\endproof

\begin{corollary}
Suppose that the estimates in the assumptions 
of the lemma hold for 
$0\leq\tau<\infty$ \ $(-\infty<\tau\leq 0)$.
Then the statement of the lemma holds
on any interval $[t_0,T]$ 
that does not contain zeros of the function $f(t)$
and $\sgn T=\sgn f(t)$, $t\in[t_0,T]$.
\end{corollary}

\section{Justification of the weak asymptotic solution}

In this section we will prove that our weak asymptotic solution is in some
sense "close" to the admissible weak solution of problem (1), (2).

The existence of the admissible weak solution in our situation
is obvious by Kruzhkov theorem (see \cite{5}, Chapter 6).

We will introduce
admissibility conditions necessary for the uniqueness of the
weak solution of considered problem.

\begin{definition} (Oleinik admissibility condition)
We say that a weak solution $u(t,x)$, $t\in \R^+$, $x\in \R$, of
problem (1), (2) is admissible if it satisfies \begin{gather*}
u_+=u(t,x^*+0)<u(t,x^*-0)=u_-.  \end{gather*} in every point of
its discontinuity.
\label{feb1114}
\end{definition}

Notice that such condition we can use
only when the function $u$ is piecevise continuous weak solutions of the
considered problem for every fixed $t\in \R^+$. In that case Definition
\ref{feb1114} is equivalent to more general Kruzhkov admissibility condition
(which can be applied on functions which are merely measurable):

\begin{definition}
We say that the weak solution $u(x,t)$, $x\in \R$, $t\in
\R^+$ of problem (1), (2) is admissible if we have \begin{gather}
\int_0^T\int_{\R}\left[ \partial_t\psi \eta(u) + \partial_x\psi
q(u)\right] dx dt+ \int_{\R}\psi(x,0)\eta(u_{0}(x))dx \geq
0,
\label{feb1214}
\end{gather} where $q(u)=\int\eta'(u)f'(u)du$ and $\eta\in
C^{1}(\R)$ is an arbitrary convex function.  \end{definition}

Using this definition, Kruzhkov proved the existance
uniqueness theorem (i.e. Theorem 6.2.2 in \cite{5}).

We will prove that weak asymptotic solution tends in $L^1$
to the admissible weak solution of problem (1), (2). By
the definition of the weak asymptotic solution
for all $\varphi\in C^\infty( [0,T];C^\infty_0(\R^1) )$
we have:
$$\int_{{\bf
R}}[u_{\eps t}+(f(u_{\eps})_x]\phi(x,t)dx=\O(\eps),$$
uniformly in $t\in [0,T]$, $T\in \R^+$.
Integrating last expression with $\int_{0}^T dt$ we
obtain:
\begin{equation}
\int_{0}^T\int_{\R}[u_{\eps
t}+(f(u_{\eps}))_x]\phi(x,t)dxdt=\O(\eps).
\end{equation}
Now letting $\eps\to 0$ we see
that $u(x,t)=w-\lim\limits_{\eps \to 0}u_{\eps}(x,t)$ is the weak solution of
(1), (2). From the construction we see that $u$ satisfy Oleinik
admissibility condition (since $u$ is obviously piecevise continuous) and
this implies Kruzhkov admissibility condition.  Furthermore, it is easy to
see that we have:  \begin{multline} \int_0^T\int_{{\R}}\left[ \partial_t\psi
\eta(u_{\eps}) + \partial_x\psi q(u_{\eps})\right] dx dt+
\int_{{\R}}\psi(x,0)\eta(u_{0}(x))dx \geq  \eps\O(1),
\label{feb1314}
\end{multline} where $q(u)=\int\eta'(u)f'(u)du$ and $\eta\in
C^{1}({\R})$ is an arbitrary convex function.

Relation~(50) holds by (49) and the smoothness of the function
$u_\varepsilon(x,t)$ for $\varepsilon>0$.

Now we can repeat the
procedure from \cite{5}, Theorem 6.2.2, page 87., to obtain:

\begin{theorem}
Let $u_{\eps}$ and $u$ satisfy (48) and (50), respectively.
There exists $s>0$
depending only on $[u_0^0,U]$ (interval in which initial data take values)
such that for every $t\in[0,T)$ and every $r>0$ we have:
$$\|u(\cdot,t)-u_{\eps}(\cdot,t)\|_{L^1(|x|<r)}
\leq (r+st)\cdot\eps\O(1).$$
\label{feb1614}
\end{theorem}


\begin{thebibliography}{99}

\bibitem{1}
V.~G.~Danilov, G.~A.~Omel$'$yanov, and V.~M.~Shelkovich,
\textit{Weak asymptotics method and interaction of nonlinear
waves}.
In: \textit{Asymptotic Methods for Wave and Quantum Problems},
M.~V.~Karasev, ed., AMS Transl., Ser.~2, Vol.~208,
33--164.

\bibitem{2}
V.~G.~Danilov,
\textit{Generalized solutions describing singularity interaction},
Int. J. of Math. and Math. Sci. \textbf{29} (2002), no.~8,
481--494.

\bibitem{3}
O.~A.~Oleinik,
\textit{On uniqueness and stability of the generalized solution
to the Cauchy problem for a quasilinear equation},
Uspekhi Mat. Nauk \textbf{14} (1959), no.~2, 165--170;
English transl. Amer. Math. Soc. Transl. (2) \textbf{23}
(1963), 285--290.

\bibitem{4}
S.~N.~Kruzhkov
\textit{First-order quasilinear equations in several
independent variables},
Math. USSR Sb. \textbf{10} (1970), 127--243.

\bibitem{5}
C.~M.~Dafermos,
\textit{Hyperbolic Conservation Laws in Continuum Physics},
Springer-Verlag, Berlin--New York, 2000.

\bibitem{6}
A.~M.~Il$'$in.
\textit{Matching of Asymptotic Expansions of Solutions
of Boundary Value Problems},
Nauka, Moscow, 1989;
English transl., Amer. Math. Soc., Providence, RI, 1992.

\bibitem{7}
V.~G.~Danilov and V.~M.~Shelkovich,
\textit{Propagation and interaction of shock waves of
quasilinear equation},
Nonlinear Studies \textbf{8} (2001), no.~1, 211--245.

\bibitem{8}
V.~P.~Maslov,
\textit{Operational Methods},
Nauka, Moscow, 1973 (In Russian);
MIR, Moscow, 1976 (In English).

\bibitem{9}
V.~P.~Maslov,
\textit{Complex WKB-Method in Nonlinear Equations},
Nauka, Moscow, 1977 (In Russian);
\textit{The Complex WKB Method for Nonlinear Equations}. I,
Birkh\"auser, Basel--Boston--Berlin, 1994 (In English).

\end{thebibliography}
\end{document}